\documentclass[12pt,twoside]{amsart}
\usepackage{amssymb}

\nonstopmode

\numberwithin{equation}{section}
\font\tengothic=eufm10 scaled\magstep 1
\font\sevengothic=eufm7 scaled\magstep 1
\newfam\gothicfam
\textfont\gothicfam=\tengothic
\scriptfont\gothicfam=\sevengothic

\DeclareMathOperator{\reg}{reg}

\DeclareMathOperator{\pnt}{\raise 0.5mm \hbox{\large\bf.}}

\DeclareMathOperator{\hhn}{\mathbb H \hbox{\rm ilb}^H(\mathbb P^n)}
\DeclareMathOperator{\hhtwo}{\mathbb H \hbox{\rm ilb}^H(\mathbb P^2)}

\DeclareMathOperator{\hredn}{\mathbb H \hbox{\rm red}^H(\mathbb P^n)}
\DeclareMathOperator{\hredthree}{\mathbb H \hbox{\rm red}^H(\mathbb P^3)}

\newtheorem{theorem}{Theorem}[section]

\newtheorem{corollary}[theorem]{Corollary}

\theoremstyle{definition}
\newtheorem{definition}[theorem]{Definition} 
\newtheorem{remark}[theorem]{Remark}
\newtheorem{example}[theorem]{Example}

\newtheorem{notation}[theorem]{Notation}

\begin{document}

\title{Families of reduced zero-dimensional schemes}
\author{Juan C.\ Migliore}
\address{Department of Mathematics, University of Notre Dame,
Notre Dame, IN 46556, USA}
\email{Juan.C.Migliore.1@nd.edu}

\thanks{Part of the work for this  paper was done while the author 
was sponsored by the National Security Agency  under Grant Number
MDA904-03-1-0071. 
}

\begin{abstract}
A great deal of recent activity has centered on the question of whether, for a
given Hilbert function, there can fail to be a unique minimum set of graded Betti
numbers, and this is closely related to the question of whether the associated
Hilbert scheme is irreducible or not.  We give a broad class of Hilbert functions
for which we show that there is no minimum, and hence that the associated Hilbert
scheme is reducible.  Furthermore, we show that the Weak Lefschetz Property holds
for the general element of one component, while it fails for every element of
another component.  
\end{abstract}

\maketitle
\tableofcontents

\section{Introduction}

In a fixed projective space $\mathbb P^n$, let $Z$ be a zero-dimensional scheme. 
We denote by $h_Z$ the Hilbert function of $Z$:
\[
h_Z(t) = \dim (R/I_Z)_t
\]
for all $t$.  We also use the analogous notation $h_{R/I_Z}(t)$ or, if $A$ is a
graded ring (e.g.\ Artinian), $h_A(t)$.

Consider the possible  zero-dimensional schemes in $\mathbb P^n$ having a fixed
Hilbert function, $H$.  Recall \cite{GMR} that the necessary and sufficient
condition for the existence of $H$ is that $H$ be a differentiable O-sequence
that stabilizes at some value $d$; i.e.\ not only $H$ but also its first
difference $\Delta H$ has to satisfy Macaulay's growth condition \cite{fsmacaulay}
(see the next section for a review), and $\Delta H (t) = 0$ for all sufficiently
large $t$.  Conversely, for any such function there always exists a reduced
zero-dimensional scheme with Hilbert function $H$.

We first review some  notation (mostly from \cite{RZ}).  We denote by
$\hhn$ the {\em postulation Hilbert scheme associated to $H$}, i.e.\ the Hilbert
scheme parameterizing zero-dimensional schemes with Hilbert function $H$. We
define $\hredn$ to be the open subscheme of $\hhn$ parameterizing the {\em
reduced} zero-dimensional schemes with Hilbert function $H$. If
$I$ is the ideal of a zero-dimensional scheme  such that
$R/I$ has Hilbert function $H$, suppose that $R/I$ has minimal free resolution
\begin{equation} \label{resol}
0 \rightarrow \bigoplus_j R^{\beta_{n,j}}(-j) \rightarrow \dots \rightarrow
\bigoplus_j R^{\beta_{1,j}}(-j) \rightarrow R \rightarrow R/I \rightarrow 0,
\end{equation}
with $\beta_{i,j} \in \mathbb N$.  We denote by $\beta^I$ the set of graded Betti
numbers of $R/I$ (or equivalently, the set of graded Betti
numbers of $I$).  We then write
\[
\mathbb B_H := \{ \beta^I : I \subset R \hbox{ and } h_{R/I} = H \},
\]
and
\[
\mathbb B'_H := \{ \beta^I : I \subset R , h_{R/I} = H \hbox{ and $I$ is
reduced} \}.
\]
$\mathbb B_H$ (resp.\ $\mathbb B'_H$) is a partially ordered set under the
inequality $\leq$, where we write $\beta^{I} \leq \beta^{I'}$ if $\beta_{i,j} \leq
\beta'_{i,j}$ (using the usual inequality for integers) for all $i,j$.

Two (related) questions may be asked, for a given Hilbert function $H$.  First, is
$\hredn$  reducible?   Second, what can we say about $\mathbb B_H$ and $\mathbb
B'_H$; and in particular, is there a unique minimum $\beta^I \in \mathbb B_H$
(resp.\ $\beta^I
\in \mathbb B'_H$) under the  partial ordering $\leq$? One can also ask similar
questions for Hilbert functions and graded Betti numbers of graded Artinian
algebras.

The two questions are related by the following theorem of Ragusa and Zappal\`a:

\begin{theorem}[\cite{RZ}] \label{RZ theorem}
If $\mathbb B_H$ (resp.\ $\mathbb B'_H$) has no unique minimum element then $\hhn$
(resp.\  $\hredn$) is reducible.
\end{theorem}

The point of this theorem is that the graded Betti numbers obey semicontinuity
for ideals in an irreducible family (cf.\ \cite{RZ}, Lemma 1.2), and within such
a family cancelation of Betti numbers can only come in consecutive terms in the
resolution (cf.\ \cite{peeva}).  So we can make the following rephrasing of
Theorem \ref{RZ theorem}: Let 
$I$ and $I'$ be ideals of  zero-dimensional schemes with the same Hilbert
function,
$H$, and assume that
$\beta^I$ and $\beta^{I'}$ are incomparable under the partial ordering.  If no
ideal $J$ exists with $\beta^J$ obtainable both from $\beta^I$ and $\beta^{I'}$ by
consecutive cancelation, then $I$ and $I'$ correspond to different components
of  $\hhn$, which is thus reducible.  If $I$ and $I'$ are both reduced then
$\hredn$ is reducible.  We will use this idea without comment in the following
sections.

\begin{definition} \label{strongly incomparable def}
We say that the Betti diagrams of $I$ and $I'$ are {\em strongly incomparable} if
the conditions of the previous paragraph are satisfied.  \qed
\end{definition}

It is known
that $\hhtwo$ is smooth and irreducible \cite{gotzmann} and that
sometimes even $\hhn$  is irreducible \cite{gotzmann2}.  On the other hand,
recently there has been a great deal of activity in this area, giving examples of
many different kinds  to show that sometimes
$\hhn$ and $\hredn$ (or related objects in the Artinian situation) are reducible,
or that $\mathbb B_H$ has no minimum element (cf.\ for instance \cite{BV},
\cite{CE},
\cite{ER}, \cite{FR}, \cite{kleppe}, \cite{RZ}, \cite{RZ2}, \cite{richert},
\cite{rodriguez}).  Some of these examples are isolated, and some (e.g.\ 
\cite{BV}, 
\cite{RZ2}, \cite{richert}) are infinite families.  While we present a new large
class of Hilbert functions for which $\hredthree$ is reducible, the overlap with
any of these previous results is minimal.

To describe the results in this paper, we first recall some notions.

\begin{definition} An Artinian graded algebra, $A$, has the {\em Weak Lefschetz
Property (WLP)} if, for a general linear form $L$, the induced multiplication
$\times L : A_i \rightarrow A_{i+1}$ has maximal rank, for all $i$.  
 A  zero-dimensional scheme is said to have WLP if its general Artinian
reduction does.  A reduced zero-dimensional scheme $Z$ has the {\em Uniform
Position Property (UPP)} if all subsets of the same cardinality have the same
Hilbert function.
\end{definition}

In the expository paper \cite{geomhilb}, the author discussed the Weak Lefschetz
Property (WLP) and the Uniform Position Property (UPP), describing both as open
conditions.  He showed in Example 3.4 that if $H$ is the Hilbert function with
first difference $(1,3,6,9,11,11,11)$ (corresponding to a set of 52 points in
$\mathbb P^3$), then $\hredthree$ is reducible.  Furthermore, in one component
no element has WLP, while in another component the general element has  WLP
(and in fact also UPP).  

The purpose of this paper is to put this example into a  much more general
framework, in the process giving a much larger class of Hilbert functions for
which $\hredthree$ is reducible.  In addition, one component has general element
with WLP, and another has no element with WLP.  We also recall (see above) that 
the Hilbert function, $H$, of a zero-dimensional scheme  is a differentiable
O-sequence.  To specify $H$, it is equivalent to specify its first difference,
$\Delta H$.  This is also known as the associated  {\em $h$-vector}.

Our main result is the following (cf.\ Corollary \ref{main corollary}):

\medskip

\noindent {\bf Theorem} {\em 
Let $H$ be a Hilbert function in four variables, for which the first difference
is of the form 
\[
\Delta H = \{ 1, 3, b_2, b_3, \dots, b_{t-1}, d, \dots, d, b_{s+2}, b_{s+3},
\dots , b_r, 0 \}.
\]
Here $t$ is simply the first degree for which the value is $=d$, and $s:= \lfloor
\frac{d-2}{2} \rfloor$.   We assume that 
$(1,3,b_2,\dots,b_{t-1},d,\dots,d)$ is again a differentiable O-sequence, and
that 
\[
\left \{ 
\begin{array}{ll}
b_{s+2} \leq d-2 & \hbox{if $d$ is even;} \\
b_{s+2} \leq d-1 & \hbox{if $d$ is odd}
\end{array}
\right.
\]
 and $b_i \geq b_{i+1}$ for all $i \geq
s+2$.  Assume furthermore that there are ``the right number" of $d$'s in the
middle. More precisely, we require that 
\[
t \leq s-1 = \left \lfloor \frac{d-3}{2} \right \rfloor.
\]
 Then $\hredthree$ is reducible.  Furthermore, on one component no element has
WLP, while on another component the general element has WLP. }

\medskip

To put this in a different context, recall the following result:

\begin{theorem}[\cite{HMNW}, Proposition 3.5]
Let $\underline{h} = (1,b_1,b_2,\dots,b_r)$ be a finite sequence of positive
integers.  Then $\underline{h}$ is the Hilbert function of a graded Artinian
$k$-algebra $R/J$ having WLP if and only if $\underline{h}$ is a unimodal
O-sequence such that the positive part of the first difference is also an
O-sequence.
\end{theorem}

Notice that the conditions on $\Delta H$ in our result above are
very similar to this necessary and sufficient condition for WLP.  The
only differences are the requirement that we have ``the right number" of $d$'s,
and that
$b_{s+2} \leq d-2$ (instead of $d-1$) when $d$ is even.   However, the condition
on the right number of $d$'s can be restrictive -- see Example \ref{ci ex}, where
we also discuss how one might extend this result.  Nevertheless, we have the nice
fact that the first difference of the Hilbert function can begin with  any
differentiable O-sequence, and  end with any non-increasing sequence that goes to
zero.  This degree of freedom is surprising: it is not difficult to construct
different-looking sets of points with the same Hilbert function, but to show that
they have strongly incomparable Betti diagrams is difficult in general.


\section{Preliminaries}

Let $R = k[x_0,\dots,x_n]$, where $k$ is an algebraically closed field.  
We first recall some basic facts about Hilbert functions, especially
those of arithmetically Cohen-Macaulay (ACM) curves in $\mathbb P^3$.  

\begin{definition} 
Let $Z \subset \mathbb P^{n-1}$ be any closed subscheme with defining (saturated)
ideal
$I = I_Z$.  We say that $Z$ is {\em arithmetically Cohen-Macaulay} (ACM) if the
coordinate ring $R/I_Z$ is a Cohen-Macaulay ring.  Note that if $Z$ is a
zero-dimensional scheme then it is automatically ACM. \qed
\end{definition}

\begin{definition} \label{ibinomexp}
The {\em $i$-binomial expansion} of the integer $a$ ($i, a >0$) is the unique
expression 
\[
a = \binom{m_i}{i} + \binom{m_{i-1}}{i-1} + \dots + \binom{m_j}{j},
\]
where $m_i > m_{i-1} > \dots > m_j \geq j \geq 1$.  \qed
\end{definition}

\begin{definition}
If $a \in {\mathbb Z}$ ($a>0$) has $i$-binomial expansion as in Definition
\ref{ibinomexp}, then we set  
\[
a^{\langle i \rangle} = \binom{m_i+1}{i+1} + \binom{m_{i-1}+1}{i} + \dots +
\binom{m_j+1}{j+1}.
\]
Note that this defines a collection of functions $^{\langle i \rangle}: \mathbb Z
\rightarrow \mathbb Z$.  \qed
\end{definition}

For example, the 5-binomial expansion of 76 is
\[
76 = \binom{8}{5}  + \binom{6}{4} + \binom{4}{3} + \binom{2}{2},
\]
so
\[ 
76^{\langle 5 \rangle} = \binom{9}{6} + \binom{7}{5} + \binom{5}{4} + 
\binom{3}{3} = 111.
\]

\begin{definition} \label{o-seq}
A sequence of non-negative integers $\{ a_i : i \geq 0 \}$ is called an {\em
O-sequence} if
\[
a_0 = 1 \ \hbox{ and } \ a_{i+1} \leq a_i^{\langle i \rangle}, 
\]
for all $i$.  An O-sequence is said to have {\em maximal growth from degree $i$ 
to degree $i+1$} if $a_{i+1} = a_i^{\langle i \rangle}$.  \qed
\end{definition}

\begin{theorem}[\cite{fsmacaulay}] \label{macaulay thm}

The following are equivalent:

\begin{itemize}
\item[(i)] $\{ a_i : i \geq 0 \}$ is an O-sequence;

\item[(ii)] $\{ a_i : i \geq 0 \}$ is the Hilbert function of a standard graded
$k$-algebra.
\end{itemize}
\end{theorem}

\begin{definition}
Given a sequence of non-negative integers $ \underline{a} = $  $\{ a_i :
i \geq 0 \}$, the {\em first difference} of this sequence is the sequence $\Delta
\underline{a} := \{ b_i \}$ defined by $b_i = a_i - a_{i-1}$ for all $i$.  (We
make the convention that $a_{-1} = 0$, so $b_0 = a_0 = 1$.)  We say that
$\underline{a}$ is a {\em differentiable O-sequence} if $\Delta \underline{a}$ is
again an O-sequence.  By taking successive first differences, we inductively
define the {\em $k$-th difference}
$\Delta^k \underline{a}$.
\qed
\end{definition}

Recall that if $H = \{ h_i : i \geq 0 \}$ is the Hilbert function of an
arithmetically Cohen-Macaulay subscheme of $\mathbb P^n$ of dimension $r$ (so the
coordinate ring $R/I$ has Krull dimension $r+1$) then the sequences $\Delta^{i}
H$, $1 \leq i \leq r+1$, are all O-sequences.  We say that $H$ is $(r+1)$-times
differentiable.

In particular, let $H$ be the Hilbert function of a non-degenerate ACM curve,
$C$, of degree $d$ in $\mathbb P^3$.    Then $\Delta H$ is the Hilbert function
of a proper hyperplane section of $C$, which is a zero-dimensional scheme of
degree $d$ in the hyperplane $\mathbb P^2$.

Another result that we will need is the following.  Let $V$ be any reduced
subscheme of $\mathbb P^n$ (not necessarily ACM) with Hilbert function $H =
(1,n+1, h_2, h_3, \dots)$.  Let $e$ be any positive integer, and define a new
sequence
$\{ e_i : i \geq 0 \}$ by $e_i = \min \{ h_i, e \}$.  It is shown in \cite{GMR}
that
$\{ e_i \}$ is the Hilbert function of some reduced set of points on $V$.  The
sequence $\{ e_i \}$ is called {\em the truncation of $H$ at $e$}.

A useful way of presenting the graded Betti numbers occurring in the minimal free
resolution of a graded module is by the so-called {\em Betti diagram}, introduced
in the computer algebra program {\tt macaulay} \cite{macaulay}, a program which
was  used to generate many examples that contributed to the theorems in this
paper.  Specifically, for a standard graded algebra $R/I$ we have the diagram
\[
\begin{array}{cccccccccccccccccccccc}
0 & : & 1 & \beta_{1,1} & \beta_{2,2} & \beta_{3,3} \\
1 & : & - & \beta_{1,2} & \beta_{2,3} & \beta_{3,4} \\
2 & : & - & \beta_{1,3} & \beta_{2,4} & \beta_{3,5} \\
3 & : & - & \beta_{1,4} & \beta_{2,5} & \beta_{3,6} \\
\vdots && \vdots & \vdots & \vdots & \vdots
\end{array}
\]
Note that the minimal free resolution has redundant, or ``ghost'' terms if and
only if two non-zero entries are adjacent diagonal (with slope 1); for instance,
if $\beta_{2,5}$ and $\beta_{3,5}$ are both non-zero, then there is a copy of
$R(-5)$ in both the second and third free modules in the resolution.


\section{Basic Construction}

For most of the remainder of this paper we will work in the ring $R =
k[x_0,x_1,x_2,x_3]$ but we will also need the ring $S = k[x_1,x_2,x_3]$. The main
results of this paper are all built on the construction that we will describe in
this section.  There are certain invariants that we will need.

\begin{notation} \label{def of t}
Consider a twice differentiable O-sequence, $H = \{ h_i : i \geq 0 \}$ and assume
that  $h_1 = 4$ and that  $\Delta H$ levels off at the value $d$.  So $ H$
is the Hilbert function of an ACM curve, $C$,  in $\mathbb P^3$ and $\Delta H$ is
the Hilbert function of its general hyperplane section.

Let $\Delta H = \{ b_i : i \geq 0 \}$ and let $\Delta^2 H = \{ c_i : i \geq 0
\}$.  Then we define the integer $t$ by 
\[
\Delta H = (1, 3, b_2, \dots, b_{t-1}, b_t = d, d, \dots)
\]
where $b_2 < b_3 < \dots < b_{t-1} < b_t = d$, 
and we note that 
\[
\Delta^2 H = (1,2, c_2,\dots,c_{t-1}, c_t, 0)
\]
with $c_t > 0$ and $\sum_{i=0}^t c_i = d$. 
The integer
$t$ has been given many different names in the literature, and we have $\reg(I_C)
= t+1$.  The finite sequence $\Delta^2 H$ is also known as the {\em $h$-vector} of
$C$. 

Now we set 
\[
s :=  \left \lfloor \frac{d-1}{2} \right \rfloor
\]
and we let $CI$ be the Hilbert function of a complete intersection of type
$(2,s)$ in four variables.  That is,
\[
\Delta CI = (1,3,5, \dots, 2s-3, 2s-1, 2s) \ \hbox{ and } \ \Delta^2 CI =
(1,2,2,\dots,2,2,1)
\]
and note that $\Delta^2 CI$ ends in degree $s$.

We make the following key assumptions:
\begin{equation} \label{key assn}
d > 3 \ \hbox{ and } \ t \leq s-1 = \left  \lfloor \frac{d-3}{2} \right \rfloor.
\end{equation}
 \qed
\end{notation}

We will assume that $d$ is even.  The case $d$ odd is entirely analogous and
will be left to the reader (but we make comments from time to time about this
case).  Note that if $d$ is even then we have  $2s = d-2$. 

\medskip

\noindent {\bf Claim 1:} $b_2 = 6$.  

Suppose otherwise.  Recall that $H$ is twice differentiable, i.e.\ $\Delta^2 H$
is an O-sequence.  We have:
\begin{itemize}
\item $b_2 \geq 3$ since $\Delta H$ is  differentiable;

\item if $b_2 = 3$ then $d = 3$ (again using differentiablility), contradicting
our assumption;

\item if $b_2 = 4$ then $\Delta^2 H$ is the sequence
\[
1, 2, 1, \dots
\]
and by Macaulay's growth condition each subsequent entry must be 1 until degree
$t$, so
$d = t+2$, which clearly violates (\ref{key assn});

\item if $b_2 = 5$ then $\Delta^2 H$ is the sequence
\[
1,2,2,c_3,\dots,
\]
where each $c_i$ is $\leq 2$ and $c_{t+1} = 0$.  Hence $d \leq 2t+1$, so
combining with (\ref{key assn}), we get
\[
\frac{d-1}{2} \leq t \leq \left \lfloor \frac{d-3}{2} \right \rfloor,
\]
a contradiction.
\end{itemize}

\noindent Therefore we have finished Claim 1.

As noted, $CI$ is the Hilbert function of a complete intersection in $R$ of type
$(2, \frac{d-2}{2})$, $\Delta^2 CI$ ends in degree $\frac{d-2}{2} =
s$, and $2s = d-2$.  Consider the following table, where the last line is obtained
by taking the difference $\Delta H - \Delta CI$. 

 {\scriptsize 
\begin{equation} \label{table of calculations}
\begin{array}{l|ccccccccccccccccccccccccccccccccccc}
\deg & 0 & 1 & 2 & 3 & 4 & \dots & t-1 & t & t+1 & \dots  & s-2 & s-1 & s &
s+1 & \dots \\
\hline
\Delta H & 1 & 3 & 6 & b_3 & b_4 & \dots & b_{t-1} & d & d & \dots & d & d & d
& d & \dots \\
\Delta CI & 1 & 3 & 5 & 7 & 9 & \dots & 2t-1 & 2t+1 & 2t+3 & \dots  & 2s-3 &
2s-1 & 2s &  2s & \dots \\ \cline{2-16}
&&& 1 & e_1 & e_2 & \dots & e_{t-3} & e_{t-2} & e_{t-1} & \dots &  5 & 3 & 2 & 2 &
\dots
\end{array}
\end{equation}
}

\noindent Since $t \leq s-1$, we verify that the sequence $\{ e_i : i \geq 0 \}$
 ends with the subsequence \linebreak $(\dots, 3,2,2,2,\dots)$, although the
columns between $t$ and $s-1$ could be redundant.  (When $d$ is odd, the sequence
$\{ e_i : i \geq 0 \}$ ends with $(\dots, 2,1,1,1,\dots)$.)

\medskip

\noindent {\bf Claim 2:} {\em $\{ e_i : i \geq 0 \}$ is an O-sequence.  }

Note that
$e_1
\leq 3$.   We want to show that for any $i$, the growth from $e_i$ to $e_{i+1}$
does not violate Macaulay's theorem.  

\medskip

\noindent {\underline{Case 1}}: $i \geq t-2$.

Then $i+2 \geq t$ so $b_{i+2} = d$ and it follows that $e_i \geq e_{i+1}$ (and in
fact $e_i > e_{i+1}$ for $t-2 \leq i \leq s-3$).  Hence $e_i$ is an O-sequence in
this range.

\medskip

\noindent {\underline{Case 2}}:  $i \leq t-3$.

Note that under
this assumption we have 
\begin{equation} \label{biei}
b_{i+2} = e_i + 2(i+2) +1 = e_i + 2i + 5
\end{equation} and  
\begin{equation} \label{binom coeff}
b_{i+2} =
\binom{i+4}{2} \hbox{ if and only if } e_i = \binom{i+2}{2}.
\end{equation}

Now note that 
\begin{equation} \label{fact1}
b_{i+3} - b_{i+2} = e_{i+1} - e_i +2.
\end{equation}
 What is the $i$-binomial expansion of $e_i$?  We have 
\[
e_i = \binom{m_i}{i} + \binom{m_{i-1}}{i-1} + \dots + \binom{m_j}{j}
\]
for some $m_i > m_{i-1} > \dots > m_j \geq j \geq 1.$  Since $\binom{a}{b} =
\binom{a}{a-b}$ for $a>b$, and taking into account (\ref{binom coeff}), we see
that without loss of generality we may write
\[
\begin{array}{rcl}
e_i & = & \displaystyle \left [ \binom{i+1}{1} + \binom{i}{1} +  \dots +
\binom{k}{1} \right ] + \left [ \binom{k-1}{0} + \dots + \binom{k-l}{0} \right ]
\\ \\
& = & [(i+1) + i +  \dots + k ] + l \\
\end{array}
\]
for some $k,l$.  Then from (\ref{biei}) we get
\[
b_{i+2} = [(i+3) + (i+2) + (i+1) + i + \dots + k] + l.
\]
Similarly we obtain
\[
\begin{array}{rcl}
e_{i+1} & = & [(i+2) + \dots + k'] + l' \\ \\
b_{i+3} & = & [(i+4) + (i+3) + (i+2) + \dots + k'] + l'.
\end{array}
\]
Now, we know that $\Delta H$ is an O-sequence, so the growth from $b_{i+2}$ to
$b_{i+3}$ satisfies Macaulay's bound.  That is,
\[
\begin{array}{rcl}
(i+4) + (i+3) + \dots + k' + l' & \leq & b_{i+2}^{\langle i+2 \rangle} \\
& = & [(i+4) + (i+3) + \dots + (k+1)] + l.
\end{array}
\]
Then
\[
\begin{array}{rcl}
e_{i+1} & = & b_{i+3} - (2i+7) \\ \\
& \leq & b_{i+2}^{\langle i+2 \rangle} - (2i+7) \\ \\
& = & e_i^{\langle i \rangle}
\end{array}
\]
proving that $\{ e_i : i \geq 0 \}$ is an O-sequence as claimed.  This complete
the proof of Claim~2.

\medskip

Now we have that all three sequences in (\ref{table of calculations}) are
infinite O-sequences.  Let $J$ be the lex-segment ideal in $S$ with Hilbert
function $\{ e_i : i \geq 0 \}$.  Note that since $J$ is lex-segment and $e_i =
2$ for all $i$ sufficiently large, we have that $J_i = \langle x_1,x_2
\rangle_{i}$ for all $i$ sufficiently large.

Let $Y$ be the reduced subscheme of $\mathbb P^3$ obtained by ``lifting" the
monomial ideal $J$ (cf. for instance \cite{MN2}).  $Y$ is reduced, the reduction
modulo $x_0$ of $I_Y$ is $J$, and the top dimensional part, $\bar Y$, of $Y$ is a
reduced curve of degree 2.  Note that the Betti diagram for $I_Y$ is the same as
the Betti diagram for $J$.  The first difference of the Hilbert function of
$Y$ is precisely $\{ e_i : i \geq 0 \}$.

The next claim is not needed for the rest of the proof, but is an interesting
observation.

\medskip

\noindent {\bf Claim 3:} {\em $\bar Y$ is ACM, and in fact it is a plane curve
consisting of two lines meeting in a point.}

(Note that when $d$ is odd, $\bar Y$ is just a line and there is nothing to prove
here.  But when $d$ is even, $\bar Y$ could in principle be two skew lines.  We
have to prove that regardless of the sequence $\{ e_i : i \geq 0 \}$, the curve
$\bar Y$ obtained in this way consists of two lines meeting in a point.)

The reduction of $I_{\bar Y}$ modulo $x_0$ is the ideal of a zero-dimensional
scheme of degree 2 (not necessarily saturated, a  priori, if $\bar Y$ is not ACM),
and this reduction contains $J$.  But in large degree
$J$ agrees with the ideal $\langle x_1,x_2 \rangle_{i}$, which is  the
saturated ideal of a zero-dimensional scheme of degree 2.  Therefore the
reduction of $I_{\bar Y}$ modulo $x_0$ agrees with $\langle x_1,x_2 \rangle$ in
all sufficiently large degrees, and hence its saturation is $\langle x_1,x_2
\rangle$.  But this ideal is supported at a point, so the curve
$\bar Y$ meets the hyperplane defined by
$x_0$ in  one point (up to multiplicity).  Hence $\bar Y$ is the union of two
lines meeting at this point, and we have proved Claim 3.

\medskip

We now choose polynomials $F \in (I_Y)_s$ and $Q \in R$ such that $Q$ is a
quadric and $(F,Q)$ is a regular sequence.  Let $V$ be the complete intersection
defined by $(F,Q)$.  Note that $h_V = CI$ and that $\deg V = d-2$.  We consider
the ideal
\[
I = Q \cdot I_Y + (F).
\]
The ideal $I$ is formed by a special case of Generalized Liaison Addition
\cite{GM4}.  The following can easily be deduced using \cite{GM4}.

\begin{enumerate}
\item $I$ is a saturated ideal defining a curve, $X$, of degree $d$.

\item The Hilbert function of $X$ is $H$.  (This comes directly from (\ref{table
of calculations}), using the methods of \cite{GM4}.)

\item $X$ is reduced, and as sets, we have $X = Y \cup V$.  Since $Q$ can be
chosen freely, $V$ avoids any 0-dimensional components of $Y$.

\item $F$ can be chosen to be a union of planes, as can $Q$, in which case $X$ is
a union of lines and reduced points.

\item \label{exact sequence} There is an exact sequence
\[
0 \rightarrow R(-s-2) \rightarrow I_Y(-2) \oplus R(-s) \rightarrow I \rightarrow
0.
\]
\end{enumerate}

Now we can say something about the Betti diagram of $I$.  
Suppose that $I_Y$ has minimal free resolution
\[
0 \rightarrow \mathbb F_3 \rightarrow \mathbb F_2 \rightarrow \mathbb F_1
\rightarrow I_Y \rightarrow 0.
\]
Notice that because of
the values of $\{ e_i : i \geq 0 \}$ in degrees $\geq s-3$ (recall the shift in
(\ref{table of calculations})), we obtain from \cite{cho-iarrobino} that 
$S/J$ (the reduction of $R/I_Y$ by $x_0$) has a non-zero socle element in degree
$s-3$ (where the ``3'' is), and in fact the socle in this degree is
1-dimensional.  Hence  since $I_Y$ and
$J$ have the same Betti diagram, we get 
\[
\mathbb F_3 = R(-s) \oplus A
\]
for some free module $A$.  Because the Hilbert function of $S/J$, $(\dots, 3, 2, 2
\dots)$ (with the 3 in degree $s-3$) has maximal growth from degree $s-2$ to
degree $s-1$, the ideal $J_{\leq s-2}$ generated by the components of degrees
$\leq s-2$ has regularity $s-2$.  In fact, since the value of the Hilbert
function of $S/J$ is 2 from that point on, there are no further generators and
without loss of generality we may substitute $J$ for $J_{\leq s-2}$.  Hence we
also have that $A$ has no summands $R(-s-1)$ or higher, and that $\mathbb F_2$ has
no summand $R(-s)$ or higher.
We get the commutative diagram
\[
\begin{array}{ccccrclcccccccccccccccccccc}
&&&&& 0 \\
&&&&& \downarrow\\
&&&& [R(-s-2) \oplus A(-2)] & \oplus & 0 \\
&&&&& \downarrow \\
&& 0 && \mathbb F_2(-2) & \oplus & 0 \\
&& \downarrow &&& \downarrow \\
&& R(-s-2) & \stackrel{\alpha'}{\longrightarrow} & \mathbb F_1(-2) & \oplus &
R(-s)
\\ && \downarrow &&& \downarrow \\
0 & \rightarrow & R(-s-2) & \stackrel{\alpha}{\longrightarrow} & I_Y(-2) & \oplus
& R(-s) &
\rightarrow  & I & \rightarrow & 0 \\
&& \downarrow &&& \downarrow \\
&& 0 &&& 0
\end{array}
\]
where $\alpha = (F,Q)$.  From the associated mapping cone and the observations
above, we get the following Betti diagram for $R/I$:

\begin{equation} \label{betti diag}
\begin{array}{cccccccccccccccccccccc}
0 & : & 1 & * & * & * \\
1 & : & - & * & * & * \\
\vdots & & \\
s-2 & : & - & * & * & * \\
s-1 & : & - & * & * & 1 \\
s & : & - & - & 1 & - \\
s+1 & : & - & - & - &  -
\end{array}
\end{equation}
Note that in particular, $I$ has no minimal generators in degree $s+2$.

\begin{remark} \label{splitting}
Because of the way that the maps go in the commutative diagram, the two 1's in
the Betti diagram do not cancel (i.e.\ the corresponding terms in the free
resolution do not split).  

However, if $I_Y$ had a minimal generator in degree $s$, and if we chose $F$ such
that it were such a generator, then $\alpha'$ would map $R(-s-2)$ isomorphically
onto a summand of $\mathbb F_1(-2)$ (corresponding to this generator).  In this
case, there would be a term $R(-s-2)$ in the first free module in the free
resolution of
$I$ coming from the mapping cone, and it would split off with the copy of
$R(-s-2)$ that we have described in the second free module.  Of course $I_Y$ does
not have a minimal generator of degree $s$, but in the next section we will
modify things so that such a generator does exist.  \qed
\end{remark}

Again, the case $d$ odd
is entirely analogous, and is left to the reader.  (This time we use a complete
intersection of type $(2, \frac{d-1}{2}$.)


\section{Main results}

In this section we make minor modifications on the construction of the last
section in order to get our results.

\begin{theorem} \label{first main result}
Let $H = \{ h_i : i \geq 0 \}$ be a  twice differentiable O-sequence over $R$,
with invariants $d, t, s$ as in Notation \ref{def of t}; in particular, we
continue to assume that
\[
t \leq s-1 = \left \lfloor \frac{d-3}{2} \right \rfloor.
\]
  Let $\bar H$
be the truncation at the value
$h_{s+1}$, so $\bar H$ is the Hilbert function of some set of $h_{s+1}$
points in
$\mathbb P^3$.  Then ${\mathbb B}'_{\bar H}$ does not have a unique
smallest element, and the postulation scheme  ${\mathbb H} \hbox{\rm red}^{\bar
H}(\mathbb P^3)$ is reducible.  Furthermore, on one
component,
$\mathcal H_1$, the general point corresponds to a set of points with WLP, while
on another component,
$\mathcal H_2$, no point corresponds to a set of points with WLP.  
\end{theorem}

\begin{proof}
We continue to assume that $d$ is even, leaving the case $d$ odd for the reader
(but occasionally remarking on the case of $d$ odd). Notice that the first
difference of $\bar H$ is
\[
\Delta \bar H = (1,3, 6, b_3,\dots,b_{t-1},d,\dots,d, 0)
\]
where the last $d$  occurs in  degree s+1.

Consider the sequence $\{ e_i' : i
\geq 0 \}$ defined by
\[
e_i' = \left \{
\begin{array}{rll}
e_i & \hbox{if $i \leq s-1$} \\
0 & \hbox{if $i \geq s$}.
\end{array} \right.
\]
Note that $\{e_i'\}$ ends with the subsequence $(\dots, 3,2,2,0)$.  Let $J'$ be
the lexsegment ideal in $S$ with Hilbert function $\{ e_i' : i \geq 0 \}$.  Note
that $J'$ has two minimal generators in degree $s$.  Let
$Y'$ be the reduced (but now zero-dimensional) subscheme of $\mathbb P^3$ obtained
by ``lifting'' the monomial ideal $J'$.  Note that $J'$ agrees with the ideal $J$
described above in degrees $\leq s-1$, so the same is true of $I_{Y'}$ and $I_Y$,
respectively.  The first difference of the Hilbert function of $Y'$ is $\{ e_i' :
i \geq 0 \}$. Since the ideal
$I_Y$ has no generator in degree
$s$,
$I_{Y'}$ must have exactly two minimal generators in degree $s$.  (When $d$ is
odd,
$I_{Y'}$ has one minimal generator in degree $s$.)  

Now we choose $F \in (I_{Y'})_s$ and $Q \in
R_2$ as above, but we choose $F$ to be a minimal generator of $I_{Y'}$.  
We  form the ideal $I' = QI_{Y'} + (F)$.  This is again a saturated ideal,
defining a reduced subscheme, $X'$, of $\mathbb P^3$ which consists of the
disjoint union of the curve defined by $(F,Q)$ (which has degree $d-2$) and $Y'$. 
Let $H'$ be the Hilbert function of $X'$.  We make the following Hilbert function
calculation using (\ref{table of calculations}):

 {\scriptsize 
\begin{equation} \label{table of calculations 2}
\begin{array}{l|ccccccccccccccccccccccccccccccccccc}
\deg & 0 & 1 & 2 & 3 & 4 & \dots & t-1 & t & t+1 & \dots  & s-1 & s &
s+1 & s+2 & s+3 & \dots \\
\hline
e_{i-2} &&& 1 & e_1 & e_2 & \dots & e_{t-3} & e_{t-2} & e_{t-1} & \dots & 3 &
2 & 2 & 0 & 0\\
\Delta CI & 1 & 3 & 5 & 7 & 9 & \dots & 2t-1 & 2t+1 & 2t+3 & \dots  &
2s-1 & 2s &  2s & 2s & 2s &  \dots \\ \cline{1-17}
\Delta H' & 1 & 3 & 6 & b_3 & b_4 & \dots & b_{t-1} & d & d & \dots  & d & d
& d & d-2 & d-2 &  \dots 
\end{array}
\end{equation}
}
Note that $I' = I_{X'}$ agrees with $I_X$ in degrees $\leq s+1$, and its Hilbert
function is $H$ up to and including degree $s+1$.  Note also that $I_{X'}$ has
only one minimal generator in degree $s+2$ (thanks to Remark \ref{splitting}),
rather than the two that one might guess by looking at (\ref{table of calculations
2}).

Thanks to
\cite{GMR}, there exists a set of points, $Z'$, on $X'$ whose Hilbert function is
the truncation of $h_{X'}$ at level $h_{s+1}$, i.e. whose Hilbert function is
$\bar H$.  The table (\ref{table of calculations 2}) shows that $I_{Z'}$ has
$d-2$ more generators in degree $s+2$ than does $I_{X'}$, so in fact $I_{Z'}$ has
$d-1$ minimal generators in degree $s+2$.

Thanks to the Betti diagram (\ref{betti diag}) and Remark
\ref{splitting}, we have the following Betti diagram for $R/I_{Z'}$:
\begin{equation} \label{betti diag for Z'}
\begin{array}{cccccccccccccccccccccc}
0 & : & 1 & * & * & * \\
1 & : & - & * & * & * \\
\vdots & & \\
s-2 & : & - & * & * & * \\
s-1 & : & - & * & * & 1 \\
s & : & - & - & - & - \\
s+1 & : & - & d-1 & * &  *
\end{array}
\end{equation}

Now we let $C$ be an ACM curve with Hilbert function $H$.  Since $\reg(I_C) =
t+1$, $R/I_C$ has Betti diagram
\begin{equation} \label{betti diag of C}
\begin{array}{cccccccccccccccccccccc}
0 & : & 1 & * & * & - \\
1 & : & - & * & * & - \\
\vdots & & \\
t-1 & : & - & * & * & - \\
t & : & - & * & * & - \\
t+1 & : & - & - & - &  -
\end{array}
\end{equation}
Let $Z$ be a finite subset of $C$ whose Hilbert function is the truncation of $H$
at level $h_{s+1}$, i.e. $h_Z = \bar H$; note that as before $I_Z$ agrees with
$I_C$ in degrees
$\leq s+1$.  Since $t \leq s-1$, the Betti diagram of $R/I_Z$ is
\begin{equation} \label{betti diag for Z}
\begin{array}{cccccccccccccccccccccccccc}
0 & : & 1 & * & * & - \\
1 & : & - & * & * & - \\
\vdots & & \\
t-1 & : & - & * & * & - \\
t & : & - & * & * & - \\
t+1 & : & - & - & - & - \\
\vdots \\
s-1 & : & - & - & - & - \\
s & : & - & - & - & - \\
s+1 & : & - & d & * &  *
\end{array}
\end{equation}
Comparing the Betti diagrams (\ref{betti diag for Z'}) and (\ref{betti diag for
Z}), one sees that $\beta'_{3,s+2} > \beta_{3,s+2}$ and $\beta_{1,s+2} >
\beta'_{1,s+2}$.  In fact it is clear that the diagrams are strongly incomparable
(in the sense of Definition \ref{strongly incomparable def}), since there can be
no element smaller than both.  The non-existence of a unique smallest element and
the reducibility of the postulation scheme then follow immediately.  $Z$
corresponds to a point in
$\mathcal H_1$ and $Z'$ corresponds to a point in $\mathcal H_2$.

The assertion about WLP comes from the observation (from the Betti diagram) that
the Artinian reduction of
$R/I_{Z'}$ has a socle element in degree $s-1$.  But the value of the Hilbert
function of this Artinian reduction in both degrees $s-1$ and $s$ is $d$, so the
map induced by a general linear form can be neither injective nor surjective
(since it has a one-dimensional kernel).  Hence this Artinian reduction does not
have WLP.  On the other hand, the Artinian reduction of $R/I_Z$ does have WLP --
 it follows from the Cohen-Macaulay property of $R/I_C$ and the fact that
$R/I_C$ agrees with $R/I_Z$ in all degrees $\leq s+1$. 

Again, the case $d$ odd is almost identical and is left to the reader.
\end{proof}

\begin{corollary}
Under the assumptions of Theorem \ref{first main result}, if
$H$ is of {\em decreasing type} (i.e.  $\Delta^2 H$ is
strictly decreasing once it starts to decrease, so $H$ is the Hilbert function of
an {\em irreducible} ACM curve in $\mathbb P^3$) then the general element of
$\mathcal H_1$ also satisfies UPP.
\end{corollary}

\begin{proof}
It is clear from the proof of Theorem \ref{first main result}.  We simply choose
$C$ to be an irreducible ACM curve and $Z$ to be a general set of points on $C$
of the right cardinality.  We do not know if the general element of $\mathcal
H_2$ has UPP, although we believe that it does not.
\end{proof}

\begin{example}
Let 
\[
\Delta H = (1, 3, 6, 10, 14, 16, 17, 17 , \dots).
\]
Then  $t = 6$, $d=17$ and $ s = 8$,  so Theorem \ref{first main result} considers
the Hilbert function
\[
\Delta \bar H = (1, 3, 6, 10, 14, 16, 17, 17, 17, 17,0).
\]
That is, 
\[
\bar H = (1, 4, 10, 20, 34, 50, 67, 84, 101, 118, 118, \dots).
\]
The corresponding postulation scheme $\mathbb H \hbox{ilb}^{\bar H}(\mathbb P^n)$,
parameterizing sets of 118 points in
$\mathbb P^3$ with Hilbert function $\bar H$, is reducible.  \qed
\end{example}

\begin{corollary} \label{main corollary}
Let $H = \{ h_i : i \geq 0 \}$ be a twice differentiable O-sequence, with
invariants $d,t,s$ as in Notation \ref{def of t}; in particular, we assume that 
\[
t \leq s-1 = \left \lfloor \frac{d-3}{2} \right \rfloor.
\]
Write $\Delta H = \{ 1, 3, b_2, b_3, \dots, b_{t-1}, d, d, \dots \}.$
Let $\bar H$ be the Hilbert function whose first difference is 
\[
\Delta \bar H = \{ 1, 3, b_2, b_3, \dots, b_{t-1}, d, \dots, d, b_{s+2}, b_{s+3},
\dots , b_r, 0 \}
\]
where 
\[
\left \{ 
\begin{array}{ll}
b_{s+2} \leq d-2 & \hbox{if $d$ is even} \\
b_{s+2} \leq d-1 & \hbox{if $d$ is odd}
\end{array}
\right.
\]
 and $b_i \geq b_{i+1}$ for all $i \geq
s+2$.  Then all the conclusions of Theorem \ref{first main result} continue to
hold for $\bar H$.
\end{corollary}

\begin{proof}
As before, we prove the case where $d$ is even, leaving to the reader the case of
$d$ odd. Our approach here is similar to that taken in \cite{HMNW}, Proposition
3.5.

We return to the monomial ideal $J'$ and the  reduced subscheme
$X'$ obtained at the beginning of the proof of Theorem \ref{first main result}. 
By the mechanism of liftings of monomial ideals (cf.\ \cite{MN2}), we may assume
that $F \in (I_{Y'})_s$ is a union of planes.  Clearly $Q$ can also be chosen to
be a (general) union of planes.  Hence $X'$ is the disjoint union of $2s = d-2$
lines (defined by the complete intersection of $F$ and $Q$) and a finite set of
points ($Y'$) whose regularity is $s$.  We saw that the Hilbert function of $X'$
has first difference
\begin{equation} \label{delta h of X'}
(1,3,6,b_3,\dots,b_{t-1}, d,\dots,d, d-2, d-2,\dots)
\end{equation}
with the last $d$ occurring in degree $s+1$.  

Note that the union, $A$, of $2s = d-2$ lines is a complete intersection of type
$(2,s)$; in particular it is ACM.  Its Hilbert function has first difference
\[
(1,3,5, \dots, d-5, d-3, d-2, d-2 \dots )
\]
where the first $d-2$ occurs in degree $s$.  
It is easy to check (e.g.\ by considering liftings) that we may order the lines
$\lambda_1,\lambda_2,\dots , \lambda_{d-2}$ such that for each $i$, $A_i
:= \lambda_1
\cup \dots \cup \lambda_i$ is ACM.  Notice that the Hilbert function of $A_i$ has
first difference which reaches value $i$ in degree $< s$, except for $A_{d-2}$,
which reaches value $d-2$ in degree exactly $s$.

For any $i$, the base locus of the linear system $|(I_{X'})_i|$ includes all the
lines $\lambda_j$.    Furthermore, in degrees $\geq s+2$, the reduction of
$R/I_{X'}$ modulo $x_0$ agrees with the reduction of $R/(F,Q)$ modulo $x_0$.  

We will prodeed inductively.  
Consider the Hilbert function
\[
1,3,b_2,\dots,b_{t-1}, d, \dots ,d,b_{s+2}, b_{s+2}, \dots.
\]
This can be obtained by the ideal
\[
I_{X_{s+2}} := I_{Z'} \cap I_{\lambda_1} \cap \dots \cap I_{\lambda_{b_{s+2}}}
\]
i.e. by taking the ideal of $X_{s+2} := Z' \cup A_{b_{s+2}}$, where $Z'$ is the
reduced zero-dimensional scheme obtained in Theorem \ref{first main result}. 
Indeed, since $\lambda_1,\dots,\lambda_{b_{s+2}}$ are contained in the base locus
of
$I_{Z'}$ in degrees $\leq s+1$, they impose no additional conditions in those
degrees.  And in degrees $\geq s+2$ the reduction modulo $x_0$ sees only the
$b_{s+2}$ lines.  We choose a finite set of points $Z_{s+2}$ whose Hilbert
function is the truncation, i.e. whose Hilbert function has first difference
\[
1,3,b_2,\dots,b_{t-1}, d, \dots ,d,b_{s+2}, 0.
\]

For the next step we proceed as we have done here, but substituting $X_{s+2}$ for
$X'$ and $Z_{s+2}$ for $Z'$, and we produce in the end a finite set of points
with Hilbert function having first difference
\[
1,3,b_2,\dots,b_{t-1}, d, \dots ,d,b_{s+2}, b_{s+3}, 0.
\]
After a finite number of steps, we are finished.  Let us denote by $W'$ this
reduced zero-dimensional scheme with the desired Hilbert function.

For the scheme $W$ (analogous to $Z$ in Theorem \ref{first main result}) having
the same Hilbert function as $W'$ but strongly incomparable Betti diagram, we
first note that a curve $C$ can be constructed with Hilbert function having first
difference
$\{ 1, 3, b_2, \dots, b_{t-1}, d, d, \dots \}$ and such that $C$ is ACM, and as
before, it is a union of lines $\lambda_1 \cup \lambda_2 \cup \dots \cup
\lambda_d$ such that $\lambda_1 \cup \dots \cup \lambda_i$ is ACM for all $i$. 
We then proceed exactly as we did above, considering truncations and slowly
removing lines, to produce a reduced zero-dimensional scheme with the desired
Hilbert function.

Now we have to check that these zero-dimensional schemes $W$ and $W'$ have
strongly  incomparable Betti diagrams.  Note that the Betti diagram for $R/I_C$ is
the same as that in (\ref{betti diag of C}), by construction.

Recall that in Theorem \ref{first main result}, we
produced zero-dimensional schemes $Z$ and $Z'$ with Hilbert function having first
difference
\[
1,3,b_2,\dots,b_{t-1}, d, \dots ,d,0
\]
where the last $d$ occurs in degree $s+1$.  We saw that the Betti diagram for
$R/I_{Z'}$ was 
\[
\begin{array}{cccccccccccccccccccccc}
0 & : & 1 & * & * & * \\
1 & : & - & * & * & * \\
\vdots & & \\
s-2 & : & - & * & * & * \\
s-1 & : & - & * & * & 1 \\
s & : & - & - & - & - \\
s+1 & : & - & d-1 & * &  *
\end{array}
\]
while the Betti diagram for $R/I_Z$ was 
\[
\begin{array}{cccccccccccccccccccccccccc}
0 & : & 1 & * & * & - \\
1 & : & - & * & * & - \\
\vdots & & \\
t-1 & : & - & * & * & - \\
t & : & - & * & * & - \\
t+1 & : & - & - & - & - \\
\vdots \\
s & : & - & - & - & - \\
s+1 & : & - & d & * &  *
\end{array}
\]
Now instead we have constructed reduced zero-dimensional schemes $W$ and $W'$
with Hilbert function having first difference
\[
1, 3, b_2,\dots,b_{t-1}, d,\dots,d,b_{s+2}, \dots.
\]
In particular, $I_{Z'}$ (resp.\ $I_{Z}$) agrees with $I_{W'}$ (resp.\ $I_{W}$) in
all degrees $\leq s+1$.  It follows that the Betti diagram for $R/I_{Z'}$ (resp.\
$R/I_{Z}$) agrees with the Betti diagraom for $R/I_{W'}$ (resp.\ $R/I_{W}$) in
all  rows above the one labelled $s$.  The fact that the Hilbert function in
degree $s+2$ is now $b_{s+2}$ rather than 0 means that instead of having $d-1$ and
$d$ minimal generators in degree $s+2$, respectively (in the above two diagrams),
we instead have $d-1-b_{s+2}$ and $d-b_{s+2}$, respectively.  But the
incomparability of the diagrams is preserved.
\end{proof}

\begin{example}
In Example 3.4 of \cite{geomhilb} it was shown that if $\bar H$ is the Hilbert
function with first difference $(1,3,6,9,11,11,11,0)$ (corresponding to a set of
52 points in
$\mathbb P^3$), then $\mathbb
H\hbox{\rm red}^{\bar H}(\mathbb P^3)$ is reducible.  Since here $d=11, t=4$ and
$s = 5$, this is easily seen to follow immediately from Theorem \ref{first main
result}.  Furthermore, by Corollary \ref{main corollary}, also the Hilbert
function $H'$ with first difference 
\[
\begin{array}{c}
(1,3,6,9,11,11,11, 9,6,3,1) \ \hbox{ or } \\
(1,3,6,9,11,11,11, 10, 8, 8, 5, 5, 5, 4, 3,3, 1)
\end{array}
\]
 has the property that $\mathbb
H\hbox{\rm red}^{H'}(\mathbb P^3)$ is reducible.

\end{example}

\begin{example} \label{ci ex}
Richert shows that even when the Hilbert function is one of a complete
intersection, the conclusions about not having a unique minimum element of
$\mathbb B_H$ (or about the reducibility of $\hredn$) may hold.  His infinite
family of examples deals with the Hilbert function of a complete intersection of
type $(m,m+1,2m+1)$.  More generally, in
\cite{RZ2} Ragusa and Zappal\`a show that
if $H$ is the Hilbert function of a complete intersection of type $(a,b,c)$ with
$a \leq b \leq c$ and $b+3 \leq c \leq ab$, and if $(a,b,c) \neq (4,4,7)$, then
$\mathbb B_H$ does not have a unique minimum element.  Using our methods it is
also possible to obtain the Hilbert function of a complete intersection, but the
``flat part" of the Hilbert function becomes rather large.  For example, taking
$a=b=7$ requires
$c=26$, and we obtain the Hilbert function
{\small 
\[
\begin{array}{lcccc}
(1,3,6,10,15,21,
28,34,39,43,46,48,49,49,49,49,49,49,49,49,49,49, 
49,49,49,49, 48, 46, 43, \\ \\  39, 34,
28, 21, 15, 10, 6, 3, 1).
\end{array}
\]}

One might wonder if we could improve the method using a different kind of
complete intersection $CI$ in our construction (recall that we used one of type
$(2,
\lfloor \frac{d-1}{2} \rfloor)$).  The idea would be to try to diminish the
number of copies of $d$ in the first difference of the Hilbert function.  If our
goal were some restricted class of Hilbert functions, such as those of complete
intersections of type
$(a,b,c)$, experimental evidence suggests that it might be possible.  But since
this case is known \cite{RZ2}, we have to explore the general case.  

Suppose, for instance, that we wanted to show that the Hilbert function with first
difference
\[
\Delta H = (1,3,6,10, 15, 19, 23, 26, 27, 28, 29, 29, \dots, 0)
\]
corresponds to a postulation scheme $\hredthree$ which is
reducible.  Using the method of this paper, we use a complete intersection of type
$(2, 14)$ and study the table

 {\small 
\begin{equation} \label{ci2}
\begin{array}{l|ccccccccccccccccccccccccccccccccccc}
\deg & 0 & 1 & 2 & 3 & 4 & 5 & 6 & 7 & 8 & 9 & 10 & 11 & 12 & 13 & 14 & 15 & \dots
\\
\hline
\Delta H & 1 & 3 & 6 & 10 &15 & 19 & 23 & 26 & 27 & 28 & 29 & 29 & 29 & 29 &
29 & 29 & \dots
\\
\Delta CI & 1 & 3 & 5 & 7 & 9 & 11 & 13 & 15 & 17 & 19 & 21 & 23 & 25 & 27 & 28 &
28 & 
\dots
\\
\cline{2-18} &&& 1 & 3 & 6 & 8 & 10 & 11 & 10 & 9 & 8 & 6 & 4 & 2 & 1 & 1 & 
\dots
\end{array}
\end{equation}
}
Here $t = 10$ and $s = 14$.  Note that the bottom row is indeed an O-sequence. 
Theorem \ref{first main result} then draws a conclusion about the Hilbert
function with first difference 
\[
(1,3,6,10,15,
19,23,26,27,28,29,29,29,29,29,29,0).
\]

But if we tried to reduce
the number of 29's by choosing, for instance, a complete intersection of type
$(4,7)$, we obtain the diagram 
{\small 
\[
\begin{array}{l|ccccccccccccccccccccccccccccccccccc}
\deg & 0 & 1 & 2 & 3 & 4 & 5 & 6 & 7 & 8 & 9 & 10 & 11 &   \dots
\\
\hline
\Delta H & 1 & 3 & 6 & 10 &15 & 19 & 23 & 26 & 27 & 28 & 29 & 29 &  \dots
\\
\Delta CI & 1 & 3 & 6 & 10 & 14 & 18 & 22 & 25 & 27 & 28 & 28 & 28 & 
\dots
\\
\cline{2-14} &&&&& 1 & 1 & 1 & 1 & 0 & 1 & 1 & 1 &
\dots
\end{array}
\]
}
and the bottom line is not an O-sequence.

However, in (\ref{ci2}) if we instead use a complete
intersection of type $(3,9)$ then we obtain
 {\small 
\[
\begin{array}{l|ccccccccccccccccccccccccccccccccccc}
\deg & 0 & 1 & 2 & 3 & 4 & 5 & 6 & 7 & 8 & 9 & 10 & 11 & 12 &  \dots
\\
\hline
\Delta H & 1 & 3 & 6 & 10 &15 & 19 & 23 & 26 & 27 & 28 & 29 & 
29 & 29 & \dots
\\
\Delta CI & 1 & 3 & 6 & 9 & 12 & 15 & 18 & 21 & 24 & 26 & 27 & 27 & 27 & \dots 
\\
\cline{2-15} &&&& 1 & 3 & 4 & 5 & 5 & 3 & 2 & 2 & 2 & 2 &
\dots
\end{array}
\]
}
One might think that the method of this paper would allow us to start with the
O-sequence $(1,3,4,5,5,3,2,2,0)$ and proceed as before to obtain a result about
the Hilbert function with first difference $(1,3,6,10, 15, 19, 23, 26, 27, 28,
29,0)$.  But recall that the complete intersection needs to include a minimal
generator for $I_{Y'}$, and the O-sequence $(1,3,4,5,5,3,2,2,0)$ does not allow a
minimal generator of degree 9.  So in fact we have to start with the O-sequence
$(1,3,4,5,5,3,2,2,2,0)$ and then we can indeed draw a conclusion about the Hilbert
function with first difference {\small $(1,3,6,10,15,
19,23,26,27,28,29,29,0).$ } This is not covered by Theorem \ref{first main
result}.

Another obstruction to extending this technique is that if $CI$ is the Hilbert
function of a complete intersection of type $(a,b)$, then in order for the
sequence  $\{ e_i \}$ to begin with a 1 we need  $a <
b$, and also we need that $a$ be strictly smaller than the initial degree of
$H$ (or some small variant which we need not make explicit here).  If we impose
this condition, and also the condition that
$a+b \geq t+3$ and that $ab$ is smaller than $d$ but ``not too much smaller" (so
that we have maximal growth of the first difference of the Hilbert function) then
very likely the technique could be extended, as illustrated above.  But the gain
would only be to shorten the string of $d$'s in the middle of the
$h$-vector.  Given the complicated
nature of these numerical conditions (we have not even addressed what we
would need to ensure that $\{ e_i : i \geq 0 \}$ is an O-sequence), and the fact
that already our result gives a very large new class of Hilbert functions for
which the postulation scheme is reducible, it does not seem worthwhile to pursue
this.
\qed
\end{example}



\begin{thebibliography}{999}

\bibitem{macaulay} D.\ Bayer and M.\ Stillman, Macaulay: A system for
computation in algebraic geometry and commutative algebra. Source and object
code available for Unix and Macintosh computers.  Contact the authors, or
download from ftp://math.harvard.edu via anonymous ftp.

\bibitem{BV} M.\ Brignone and G.\ Valla, {\em On the resolution of certain level
algebras},  Comm. Algebra {\bf 32} (2004), no. 11,
4221--4245.

\bibitem{CE} H.\ Charalambous and G.\ Evans, {\em Resolutions with a given
Hilbert function}, in: Commutative algebra: syzygies, multiplicities, and 
birational algebra (South Hadley, MA, 1992), 19--26,  Contemp. Math., 159, 
Amer. Math. Soc., Providence, RI, 1994. 

\bibitem{cho-iarrobino} Y.\ Cho and A.\ Iarrobino, {\em Hilbert functions and
level algebras}, J.\ Algebra {\bf 241} (2001), 745--758.

\bibitem{ER} E.G.\ Evans and B.\ Richert, {\em Possible resolutions for a given
Hilbert function},  Comm. Algebra {\bf 30} (2002), no. 2,
897--906.

\bibitem{FR} C.\ Francisco and B.\ Richert, {\em Lex-plus-powers ideals},
preprint.

\bibitem{GMR} A.V.\ Geramita, P.\ Maroscia and L.\ Roberts, {\em The Hilbert
Function of a Reduced $k$-Algebra}, J.\ London Math.\ Soc.\ {\bf 28} (1983),
443--452.

\bibitem{GM4} A.V.\ Geramita and J.\ Migliore, {\em A Generalized Liaison
Addition}, J.\ Alg.\ 163 (1994), 139--164.

\bibitem{gotzmann} G.\ Gotzmann, {\em A stratification of the Hilbert scheme of
points in the projective plane}. Math. Z. {\bf 199} (1988), no. 4, 539--547. 

\bibitem{gotzmann2} G.\ Gotzmann, {\em Some irreducible Hilbert schemes}, Math.
Z. {\bf 201} (1989), no. 1, 13--17.

\bibitem{kleppe} J.\ Kleppe, {\em Maximal families of Gorenstein algebras},
to appear in Trans.\ Amer.\ Math.\ Soc.

\bibitem{HMNW} T.\ Harima, J.\ Migliore, U.\ Nagel and J.\ Watanabe, {\em The 
Weak and Strong Lefschetz Properties for Artinian $k$-algebras}, J.\ Algebra
{\bf 262} (2003), 99-126. 

\bibitem{fsmacaulay} F.S.\ Macaulay, {\em Some properties of enumeration in
the theory of modular systems}, Proc.\ Lond.\ Math.\ Soc.\ {\bf 26} (1927) 
531--555.

\bibitem{geomhilb} J.\ Migliore, {\em The geometry of Hilbert functions},
preprint.

\bibitem{MN2} J.\ Migliore and U.\ Nagel, {\em Lifting Monomial Ideals}, Comm.\
Algebra {\bf 28}  (2000) (special volume in honor of Robin Hartshorne),
5679-5701. 

\bibitem{peeva} I.\ Peeva, {\em Consecutive cancellations in Betti numbers}, 
Proc. Amer. Math. Soc. {\bf 132} (2004), no. 12, 3503--3507

\bibitem{RZ} A.\ Ragusa and G.\ Zappal\`a, {\em On the reducibility
of the postulation Hilbert scheme}, to appear in Rend.\ Circ.\ Mat.\ Palermo.

\bibitem{RZ2} A.\ Ragusa and G.\ Zappal\`a, {\em Looking for minimal graded Betti
numbers}, preprint.

\bibitem{richert} B.\ Richert, {\em Smallest graded Betti numbers}, J.\
Algebra {\bf 244} (2001), 236--259.

\bibitem{rodriguez} M.\ Rodriguez, 
{\em Ideals attaining a given Hilbert function}, 
Illinois J. Math. {\bf 44} (2000), no. 4, 821--827.


\end{thebibliography}
\end{document}